\documentclass[12pt]{amsart}
\usepackage{amsmath,amssymb,amsthm}
\usepackage{psfig}

\newtheorem{thm}{Theorem}
\newtheorem{prop}{Proposition} 
\newtheorem{rem}{Remark}

\newcommand{\ga}{\gamma}

\newcommand{\rc}{\mathcal R}

\newcommand{\Z}{\mathbb Z}
\newcommand{\s}{\sigma}

\renewcommand{\a}{\alpha}
\renewcommand{\b}{\beta}

\newcommand{\C}{\mathbb C}

\renewcommand{\d}{\delta}
\newcommand{\D}{\Delta}
\newcommand{\e}{\varepsilon}

\newcommand{\V}{U_h(gl_2)}

\newcommand{\U}{\mathcal U}
\newcommand{\ue}{{\mathcal U}_{\varepsilon}}

\newcommand{\oti}{\otimes}
\newcommand{\uto}{u_2}
\newcommand{\tuto}{\tilde{u}_2}
\newcommand{\vto}{v_2}
\newcommand{\tvto}{\tilde{v}_2}
\newcommand{\uo}{u_1}
\newcommand{\tuo}{\tilde{u}_1}
\newcommand{\vo}{v_1}
\newcommand{\tvo}{\tilde{v}_1}

\begin{document}

\setlength{\baselineskip}{16pt}

\title[Braiding
for the quantum $gl_2$]{Braiding
for the quantum $gl_2$ at roots of unity}
\author{R. Kashaev}
\address{Section de math\'ematiques, Universit\'e de Gen\`eve
CP 240, 2-4 rue du Li\`evre, CH 1211 Gen\`eve 24, Suisse}
\email{Rinat.Kashaev@math.unige.ch}
\author{N. Reshetikhin}
\address{Department of mathematics, University of California,
Berkeley, CA 94720, USA} \email{reshetik@math.berkeley.edu}

\begin{abstract}In the preceding papers \cite{KR1,KR2} we started 
considering the categories of tangles with flat $G$-connections 
in their complements, where $G$ is a simple complex algebraic group. 
The braiding (or the commutativity
constraint) in such categories satisfies the holonomy Yang--Baxter
equation and it is this property which is essential for
our construction of invariants of tangles with flat $G$-connections 
in their complements. In this paper, to any pair of irreducible
modules over the quantized universal enveloping algebra of $gl_2$ at a
root of unity, we associate  a solution of the holonomy Yang--Baxter
equation.
\end{abstract}
\maketitle
\tableofcontents

\section*{Introduction}

The representation theory of quantum groups at roots of unity
which we will use was developed in \cite{DCK} and \cite{DCKP}.
These quantum groups are finite dimensional over the central Hopf
subalgebra generated by $\ell$-th powers of generators and root
elements. This central Hopf subalgebra is usually denoted by
$Z_0$.

Quansitriangular Hopf algebras and, more generally, tensor
categories can be used for construction of invariants of links and
3-manifolds (see \cite{T1} and references therein). In a similar
way tensor categories fibered over a braided group can be used for
construction of invariants of tangles with flat connections in their
complements \cite{T2,KR2}. Less functorially one can say that
such invariants can be constructed from any system of solutions to
the holonomy Yang--Baxter equation \cite{KR1}.

Let $\ue$ be the quantized universal enveloping algebra of $gl_2$
at a root of unity $\e$ of degree $\ell$. The first step towards
understanding the braiding structure for quantized universal
enveloping algebras was done in \cite{R,Ga}. In this paper
we will show that every pair of irreducible modules over $\ue$
defines a solution to the holonomy Yang--Baxter equation for the
Lie group $GL_2^*$. The braiding on this group is birational and
it is given by the factorization mapping.

Moreover, one can show that the category of modules over $\ue$,
where objects are finite dimensional modules irreducible over
$Z_0$, is a tensor category fibered over $GL_2^*$. One can also
show that this is true for any factorizable Lie algebra related to a
simple Lie algebra: the category of finite dimensional models
irreducible over $Z_0$ is a tensor category fibered over $Z_0$.
These topics will be the subject of a separate publication.

N. R. would like to thank Claudio Procesi for many illuminating
discussions. The work of N.R. was partly supported by the NSF grant 
DMS-0070931
and by the Alexander von Humbolt foundation. Both authors were
partly supported by the Swiss National Foundation.

\section{Quantized Universal enveloping algebra of
$gl_2$}

\subsection{ The algebra $\V$}

\subsubsection{}The algebra $\V$ over the ring $\C[[h]]$ is generated
by elements $H,G,X$, and $Y$ with the defining relations
\begin{equation}\label{u1}
[H,G]=0, [H,X]=2X,  [H,Y]=-2Y,
\end{equation}
\begin{equation}\label{u2}
[G,X ]=2X, [G,Y]=-2Y,
\end{equation}
\begin{equation}\label{u3}
 [X,Y]=\frac{e^{\frac{hH}{2}}-e^{-\frac{hG}{2}}}{e^{\frac{h}{2}}-e^{-\frac
 {h}{2}}}
\end{equation}

The Hopf algebra structure on $\V$ is defined by the action of the
comultiplication on the generators:
\begin{eqnarray}
\D H&=&H\otimes 1+1\otimes H, \ \D G=G\otimes 1+ 1\otimes G, \\
\D X&=&X\otimes e^{\frac{hH}{2}}\otimes 1 + 1\otimes X, \ \D Y= Y
\otimes 1 +e^{-\frac{hG}{2}}\otimes Y
\end{eqnarray}
Elements $H, G, X$ and $Y$  are deformed versions of $2e_{11}, \
2e_{22}, \ e_{12}$ and $e_{21}$ in $gl_2$.

The algebra $\V[h^{-1}]$ is isomorphic to the Drinfeld double of the quantized
universal enveloping algebra $U_h(b)[h^{-1}]\subset
U_h(sl_2)[h^{-1}]$ (see Appendix~\ref{double}). As the double of a Hopf
algebra it is quasitriangular with the universal $R$-matrix
\cite{Dr} given by the canonical element in
$U_h(b)\widehat{\otimes}U_h(b)^\vee$ where $U_h(b)^\vee$ is a dual
counterpart of $U_h(b)$ and $\widehat{\otimes}$ is the tensor
product completed over formal power series in $h$:
\begin{equation}\label{Rh}
R=\exp(\frac{h}{4} H\otimes G) \prod_{n\geq 0}
(1-e^{-\frac{h}{2}}(e^{\frac{h}{2}}-e^{-\frac {h}{2}})^2X\otimes Y
e^{-nh})^{-1}
\end{equation}
This is the element of $\V^{\otimes 2}[[h]]$ which one should
consider as a formal power series in $h$.

Element $R$ as the universal $R$-matrix satisfies the following
identities:
\begin{equation}\label{br}
R \D(a)R^{-1}= \sigma\cdot\D(a)
\end{equation}
\[
(\D \oti id)(R)= R_{13} R_{23}
\]
\[
(id \oti \D)(R) =R_{13} R_{12}
\]
and, in particular, it satisfies the Yang--Baxter equation
\[
R_{12}R_{13}R_{23}=R_{23}R_{13} R_{12}
\]

\subsection{The inner automorphism $\rc$}
Define the inner automorphism $\rc: \V^{\oti 2}[[h]]\to\V^{\oti
2}[[h]]$ as
\begin{equation}\label{rc}
\rc(x\oti y)=R(x\oti y) R^{-1}
\end{equation}

It is easy to compute the action of $\rc$ on the generators.

\begin{thm} The following identities hold true:
\begin{equation}\label{rc-1}
\rc(1\oti e^{\frac{hH}{2}})=(1\oti
e^{\frac{hH}{2}})(1-e^{\frac{h}{2}}(e^{\frac{h}{2}}-e^{-\frac
{h}{2}})^2e^{-\frac{hH}{2}}X\otimes Ye^{\frac{hG}{2}})^{-1}
\end{equation}
\[
\rc(1\oti e^{\frac{hG}{2}})=(1\oti
e^{\frac{hG}{2}})(1-e^{\frac{h}{2}}(e^{\frac{h}{2}}-e^{-\frac
{h}{2}})^2e^{-\frac{hH}{2}}X\otimes Ye^{\frac{hG}{2}})^{-1}
\]
\[
\rc(X\oti 1)=X\oti e^{\frac{hG}{2}}
\]
\[
\rc(1\oti Y)=e^{-\frac{hH}{2}}\oti Y
\]
\end{thm}
The proof is immediate by the use of the commutation relations
between the generators and the equation
\[
f(zq;q)=(1-z)f(z;q)
\]
for the function $f(z;q)=(1-z)^{-1}(1-zq)^{-1}(1-zq^2)^{-1}\dots$.

The action of $\rc$ on elements $1\oti X$, $Y\oti 1$,
$e^{\frac{hH}{2}}\oti 1$, and $e^{\frac{hG}{2}}\oti 1$ can be
derived from the formulae above and from identity (\ref{br}).

The Yang--Baxter equation for $R$ implies the Yang--Baxter equation
for $\rc$:
\begin{equation}\label{rc-ybe}
\rc_{12}\cdot\rc_{13}\cdot\rc_{23}=\rc_{23}\cdot\rc_{13}\cdot\rc_{12}
\end{equation}

\section{The algebra $\U$}

\subsection{The algebra}The algebra $\U$ is generated over $\C[t,t^{-1}]$ by
elements $K,L,E$ and $F$ with the following defining relations
\[
KL=LK, \ KE=t^2EK, \ KF=t^{-2}FK,
\]
\[
LE=t^2EL, \ LF=t^{-2}FL,
\]
\[
EF-FE=(t-t^{-1})(K- L^{-1})
\]
The center of $\U$ is generated freely by Laurent polynomials in
$KL^{-1}$ and
\begin{equation}\label{c}
c=EF+Kt^{-1}+L^{-1}t
\end{equation}
This is a Hopf algebra with comultiplication
\[
\D(K)=K\oti K, \ \D(L)=L\oti L,
\]
\[
\D(E)=E\oti K+1\oti E, \ \D(F)=F\oti 1+L^{-1}\oti F \ .
\]

The map $\phi: \U\to \V$ acting on the generators as
\[
\phi(K)=\exp(\frac{hH}{2}), \ \phi(L)=\exp(\frac{hG}{2}),
\phi(t)=e^{\frac{h}{2}}
\]
\[
\phi(E)=(e^{\frac{h}{2}}-e^{-\frac {h}{2}})X, \
\phi(F)=(e^{\frac{h}{2}}-e^{-\frac {h}{2}})Y
\]
extends to a homomorphism of Hopf algebras.

\subsection{The outer automorphism $\rc$}\label{raut}The algebra $\U$ 
is not quasitriangular. Instead, there
is an outer automorphism of the division ring $\overline{\U^{\oti 2}}$
of $\U^{\oti 2}$ which we denote by the same letter $\rc$ as the
automorphism (\ref{rc}) which acts on the generators as
\[
\rc(1\oti K)=(1\oti K)(1-tK^{-1}E\oti FL)^{-1}
\]
\[
\rc(1\oti L)=(1\oti L)(1-tK^{-1}E\oti FL)^{-1}
\]
\[
\rc(E\oti 1)=E\oti L
\]
\[
\rc(1\oti F)= K^{-1}\oti F
\]
Define its action on the elements $K\oti 1$, $L\oti 1$, $K\oti 1$,
and $L\oti 1$ through the condition that
\begin{equation}\label{rcom}
\rc(\D (a))=\sigma\circ \D(a)
\end{equation}
where $a$ is any of the generators of $\U$.

It is clear that the homomorphism $\phi$ brings the outer
automorphism $\rc$ to (\ref{rc}). It also satisfies the
Yang--Baxter relation (\ref{rc-ybe}).

\section{The algebra $\ue$}
\subsection{The algebra}
Let $\e$ be a primitive root of $1$ of an odd degree $\ell$
(this is a technical assumption, convenient because then $\e^{2n}$
runs through all possible $\ell$-th roots of 1 for $n=1,\dots,
\ell$). Denote by $\ue$ the specialization of $\U$ to $t=\e$. The
following theorem was proven in \cite{DCK} for any simple Lie
algebra.

\begin{thm}\label{2}
\begin{itemize}
\item Elements $E^\ell$, \ $F^\ell$, \ $K^{\pm\ell}$, and
$L^{\pm\ell}$ are central in $\ue$. Denote by $Z_0$ the central
subalgebra in $\ue$ which they generate.

\item $Z_0$ is a Hopf subalgebra with the comultiplication
\[
\Delta(K^\ell)=K^\ell\oti K^\ell, \ \Delta(L^\ell)=L^\ell\oti
L^\ell,
\]
\[
\Delta(E^\ell)=E^\ell\oti K^\ell+1\oti E^\ell, \
\Delta(F^\ell)=F^\ell\oti 1+L^{-\ell}\oti F^\ell .
\]

\item The algebra $\ue$ is a free $Z_0$-module of dimension
$\ell^4$.

\item The center $Z(\ue)$ is generated by $Z_0$ and by the element
(\ref{c}) modulo the relation
\[
\prod_{j=0}^{\ell-1}(c-K\e^{j+1}-L^{-1}\e^{-j-1})=E^\ell F^\ell
\]

\item Let $a,b,c,d$ be coordinates on the group $B_+\times B_-$
such that for $b_\pm\in B_\pm$ we have:
\[
b_+=\left(\begin{array} {cc}1& b \\ 0& a\end{array}\right) \ , \
b_-=\left(\begin{array} {cc}d& 0 \\ c& 1\end{array}\right)
\]
Then the map $F^\ell\to b$, $E^\ell\to -c d^{-1}$, $K^\ell\to a$,
$L^\ell\to d^2$ is an isomorphism of Hopf algebras $Z_0\to
C(B_+\times B_-)$

\item $\ue$ is semisimple over a Zariski open subvariety in
$Spec(Z_0)\simeq B_+\times B_-$.

\end{itemize}
\end{thm}

\subsection{$Z_0$-irreducible quotients}Let $x\in GL_2^*$ be an irreducible
$Z_0$-character and $I_x\subset \ue$ be the corresponding ideal.
The quotient algebra
\[
A_x=\ue/I_x
\]
is finite-dimensional of dimension $\ell^4$. There are three
natural structures of a left module on $A_x$. For $a\in \ue$
denote by $[a]$ the class of $a$ in $A_x$. Then the three
actions are:
\begin{itemize}
\item $\pi(a)[b]=[ab]$, \item $\phi(a)[b]=[bS(a)]$, \item
$\psi(a)[b]=[bS^{-1}(a)]$.
\end{itemize}

Assume that $x\in GL_2^*$ is generic, i.e. that $A_x$ is
semisimple. Fix an isomorphism of algebras $\phi_x: A_x \simeq
\oplus_{i=1}^n Mat(k_i)$. For the algebra $\ue$ it is known
\cite{DCK,DCKP} that $n=\ell$ and $k_i=\ell$ for all
$i=1,\dots, \ell$. Define
\[
t: A_x\to \C, \ t(a)=\sum_{i=1}^\ell t_i Tr(\phi^i_x(a))
\]
where $Tr$ is the matrix trace in $Mat(k_i)$ and $\phi^i_x:
A_x\to Mat(k_i)$ is the $i$-th component of $\phi_x$.
It is clear that $t(a)$ does not depend on a particular choice of
$\phi_x$. Indeed, any other such isomorphism differs from $\phi_x$
by an inner automorphism of $\oplus_{i=1}^n Mat(k_i)$. Since
trace is cyclically invariant, the value of $t(a)$ for such
isomorphism will be the same as for $\phi_x$.
Thus, for generic $x$ we have an invariant bilinear form on
$A_x$:
\[
(a,b)=t(ab) \ ,
\]
It is non-degenerate if $t_i\neq 0$ for each $i=1,\dots, n$.

Fix a scalar product on $A_x$ as above. It
gives an isomorphism of vector spaces $A_x^*\simeq A_x$.
It is easy to verify that the the pairing between $\ue$-modules
$(A_x, \phi)$ and $(A_x, \pi)$ given by the map  $e_x: (A_x,
\phi)\otimes (A_x, \pi)\to \C$
\begin{equation}\label{e}
e_x: a\otimes b \to t(ab)
\end{equation}
is $\ue$-invariant with respect to the diagonal action. Indeed,
by using Sweedler's notation $\Delta(c)=\sum_c c^{(1)}\otimes
c^{(2)}$ for the action of the comultiplication on element $c$,
we have:
\[
e_x(\sum_c aS(c^{(1)}) \otimes c^{(2)}b)=\sum_c
t(aS(c^{(1)})c^{(2)}b)=\e(c)t(ab)
\]
Similarly
\[
e_x(\sum_c c^{(1)}a\otimes bS^{-1}(c^{(2)}))=\sum_c
t(c^{(1)}abS^{-1}(c^{(2)}))=\e(c)t(ab)
\]
Therefore, the map $e_x: (A_x, \pi)\otimes
(A_x, \psi)\to \C$ defined by (\ref{e}) is also $\ue$-invariant.

Consider the mapping $i_x: \C\to A_x\otimes A_x$ acting at $1$ as
\[
i_x(1)\mapsto\sum_i e_i\otimes e^i \,
\]
It is easy to see that it is a morphism of $\ue$-modules $\C\to
(A_x,\pi)\otimes (A_x,\phi)$. Indeed, let $a_i^j=t(ae_ie^j)$, then
\begin{multline}
\sum_a\sum_ia^{(1)}e_i\oti e^iS(a^{(2)})=\sum_a\sum_{i,j,k}
(a^{(1)})_i^j\otimes (S(a^{(2)})_k^ie^k=\\
\sum_a\sum_{j,k}
(a^{(1)}S(a^{(2)})_k^je_j\otimes e^k=\e(a)\sum_i e_i\oti e^i
\end{multline}
which implies the first statement. Similarly, the mapping 
$\C\to (A_x,\psi)\oti (A_x,\pi)$ is also a morphism of
$\ue$-modules.

\subsection{The action of $\rc$ on $Z_0\oti Z_0$} \begin{thm}\label{3}
The subspace $Z_0\oti Z_0\subset \ue\otimes
\ue$ is invariant with respect to the action of the automorphism $\rc$.
\end{thm}
\begin{proof}
From the action  of $\rc$ on the generators of $\ue$  and from
the relations between generators we have:
\[
\rc(1\otimes K^\ell)=(1\otimes K^\ell)(1+K^{-\ell}E^\ell\otimes
F^\ell L^\ell)^{-1}  \,
\]
\[
\rc(1\otimes L^\ell)=(1\otimes L^\ell)(1+K^{-\ell}E^\ell\otimes
F^\ell L^\ell)^{-1}  \,
\]
\[
\rc(E^\ell\otimes 1)=E^\ell\otimes L^\ell  \,
\]
\[
\rc(1 \otimes F^\ell)=K^{-\ell}\otimes F^{\ell}  \ .
\]
The comultiplication acts on $\ell$-th powers of the generators as
\[
\D(K^\ell)=K^\ell\oti K^\ell, \ \ \D(L^\ell)=L^\ell\oti L^\ell \,
\]
\[
\D(E^\ell)=E^\ell\oti K^\ell+1\oti E^\ell \,
\]
\[
\D(F^\ell)=F^\ell\oti 1+L^{-\ell}\oti F^\ell \ .
\]
These formulae and the defining property $\rc(\D(a))=\sigma\circ
\D(a)$ imply that
\begin{eqnarray}\label{COP}
\rc(K^\ell \oti K^\ell)&=&K^\ell\oti K^\ell,\\
\rc(L^\ell \oti L^\ell)&=&L^\ell\oti L^\ell,\\
\rc(E^\ell\oti K^\ell+1\oti E^\ell)&=&K^\ell\oti E^\ell+E^\ell\oti
1,\\
\label{COPP}
\rc(E^\ell\oti 1+L^{-\ell}\oti F^\ell)&=&F^\ell\oti L^{-\ell}+1
\oti F^\ell,
\end{eqnarray}
These formulas describe implicitly the action of $\rc$ on the rest 
of the generators
of $Z_0\oti Z_0$. In particular, it is clear that the image is in
$Z_0\oti Z_0$.
\end{proof}
Translating the action of $\rc$ on the generators of $Z_0\oti Z_0$ through
the identification of $Z_0$ and $C(GL_2^*)$ we have the following
statement.
\begin{thm}\label{4}The automorphism $\rc$ is the pull-back of the 
birational correspondence
$\beta: GL_2^*\times GL_2^*\to GL_2^*\times GL_2^*$ defined as
follows:
\begin{equation}\label{beta}
\beta: (x,y)\mapsto (x_L(x,y), x_R(x,y))
\end{equation}
where
$I(x_L(x,y))=x_-I(y)x_-^{-1}$ and
$I(x_R(x,y))=(x_L(x,y))_+^{-1}I(x)(x_L(x,y))_+$.
Here $I: GL_2^*\to GL_2$ is 
the factorization mapping $I(x_+,x_-)=x_+x_-^{-1}$.
\end{thm}
\begin{proof} Let $a_1,b_1,c_1,d_1$ and $a_2,b_2,c_2,d_2$ be
coordinates of points $x, y\in GL_2^*$ respectively (as in
Theorem \ref{2}. Let $\a_1,\b_1,\ga_1,\d_1$ and
$\a_2,\b_2,\ga_2,\d_2$ be coordinates of points $x_R,x_L\in
GL_2^*$. By definition, the pull-back of the mapping $\beta$
is $\rc$. From the explicit action of $\rc$ on the coordinate
functions we have
\[
\a_2=a_2(1-c_1b_1b_2d_2/a_1d_1d_2)^{-1}
\]
\[
\d_1=d_1(1-c_1b_1b_2d_2/a_1d_1d_2)^{-1}
\]
\[
\ga_1\d_1^{-1}=c_1d_1^{-1}d_2
\]
\[
\b_2=b_2a_1^{-1}
\]
\[
\a_1\a_2=a_1a_2
\]
\[
\ga_1\d_1^{-1}\a_2+\ga_2\d_2^{-1}=a_1c_1d_1^{-1}+\ga_2\d_2
\]
\[
\b_1+\d_1^{-1}=b_1d_2^{-1}+b_2
\]

Now, it remains to prove that this mapping can be written as
(\ref{beta}). This a simple linear algebra exercise.
\end{proof}

Let $I_x$ be the ideal in $\ue$ corresponding to the irreducible
$Z_0$-character $x\in GL_2^*$.  Theorems \ref{3} and \ref{4}
have an important implication. The mapping $\rc$ induces an
isomorphism of algebras
\[
\rc(x,y): A_x\oti A_y\to A_{x_R(x,y)}\oti A_{x_L(x,y)} \,
\]
This mapping is also an isomorphism of the tensor product of left
$\ue$-modules.

\section{Braiding for irreducible representations of $\ue$.}

\subsection{Restriction of the braiding to irreducible
representations} Let $Z_1$ be the central subalgebra generated by the
central elements $c$ and $KL^{-1}$. From the definition of $\rc$
we have
\[
\rc(c_1\oti c_2)=c_1\oti c_2
\]
for any $c_{1,2}\in Z_1$. For generic $x$ the algebra $A_x$ is
semisimple and its irreducible representations are separated by
eigenvalues of $c$. Let $A_x^i$ be one of the irreducible
representations of $A_x$.

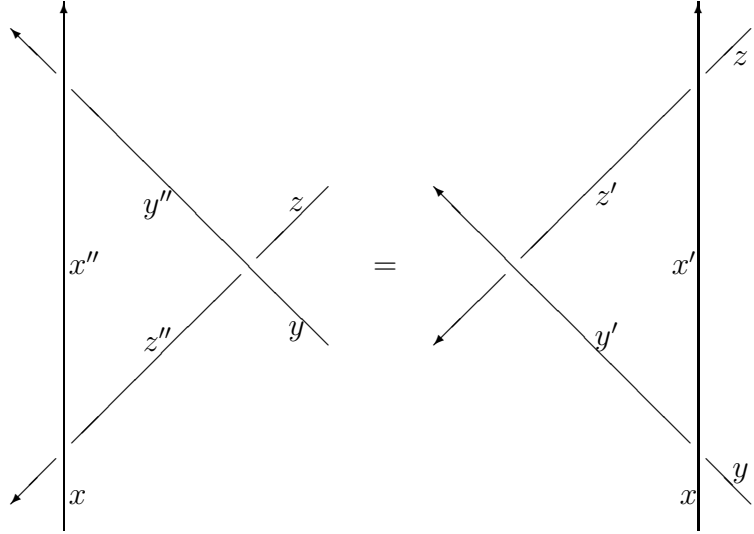
\begin{figure}
\begin{picture}(300,200)

\put(30,0){\vector(0,1){200}}
\put(130,70){\line(-1,1){97}}
\put(27,173){\vector(-1,1){17}}
\put(130,130){\line(-1,-1){27}}
\put(97,97){\line(-1,-1){64}}
\put(27,27){\vector(-1,-1){17}}
\put(270,0){\vector(0,1){200}}

\put(32,10){$x$}
\put(32,97){$x''$}
\put(115,74){$y$}
\put(60,120){$y''$}
\put(115,121){$z$}
\put(60,68){$z''$}

\put(147,97){$=$}

\put(263,10){$x$}
\put(260,97){$x'$}
\put(283,176){$z$}
\put(231,123){$z'$}
\put(283,20){$y$}
\put(231,70){$y'$}

\put(290,10){\line(-1,1){17}}
\put(267,33){\vector(-1,1){97}}
\put(290,190){\line(-1,-1){17}}
\put(267,167){\line(-1,-1){64}}
\put(197,97){\vector(-1,-1){27}}

\end{picture}
\caption{\label{hYBE} Holonomy Yang--Baxter equation.}
\end{figure}

Since $\rc$ acts trivially on $Z_1\otimes Z_1$, the isomorphism
$\rc(x,y): A_x\oti A_y\to A_{x_R(x,y)}\oti A_{x_L(x,y)}$ restricts
to the subalgebra $A_x^i\otimes A_y^j$ and induces an algebra
isomorphism
\[
\rc^{i,j}(x,y): A_x^i\otimes A_y^j\to A^i_{x_R(x,y)}\oti
A^j_{x_L(x,y)} \  .
\]
Algebras $A_x^i$ are fibers of a   bundle of algebras
with the fiber $Mat(\ell\times \ell)$. This bundle is trivializable
over sufficiently small  neighborhood of 1. Fix such trivialization, i.e.
for each $x$ in this neighborhood fix an algebra isomorphism
$\phi^i_{x}: A^i_x\simeq Mat(\ell)$. Then the mapping
$\rc^{i,j}(x,y)$ induces an automorphism of $Mat(\ell)^{\otimes
2}$. Since all automorphisms of a matrix algebra are inner, there
exists $R^{i,j}(x,y)\in Mat(\ell)^{\otimes 2}$ such that
\begin{equation}\label{inbr}
(\phi^i_{x_R(x,y)}\oti \phi^j_{x_L(x,y)})\circ \rc^{i,j}(x,y)\circ
((\phi^i_{x})^{-1}\oti (\phi^i_{y})^{-1})(A\oti B)=
R^{i,j}(x,y)(A\oti B)R^{i,j}(x,y)^{-1} \ ,
\end{equation}
for any $A, B\in Mat(\ell)$.
The Yang--Baxter relation for $\rc$ implies the projective holonomy
Yang--Baxter  equation for $R^{i,j}(x,y)$,
\begin{equation}\label{hYBEc}
R^{i,j}(x'',y'')_{12}R^{i,k}(x,z'')_{13}R^{j,k}(y,z)_{23}=c^{i,j,k}(x,y,z)
R^{j,k}(y',z')_{23}R^{i,k}(x',z)_{13}R^{i,j}(x,y)_{12}
\end{equation}
where $x, x',x'',y,y'$ etc are $GL_2^*$-colorings of the diagrams of
Fig.~\ref{hYBE} (see \cite{KR1,KR2}). In terms of the 
mapping $\beta:(x,y)\mapsto
(x_L(x,y), x_R(x,y)$  the arguments in this equation are given by
\[
z''=x_L(y,z),\ y''=x_R(y,z),\ x''=x_R(x,z''), \ x'=x_R(x,y),\  
y'=x_L(x,y),\ z'=x_L(x',z).
\]

The function $c^{i,j,k}(x,y,z)$ can be determined from the
comparison of determinants of the left and right hand sides of the
equation. We will describe a normalization of $R^{i,j}(x,y)$
for which $c^{i,j,k}(x,y,z)=1$.

\subsection{Irreducible representations of $\ue$.}
Now we fix the isomorphisms $\phi^i_x: A^i_x\to Mat(\ell)$ by
choosing a specific realization of irreducible representations of
$\ue$ . This realization is known as a cyclic representation.

Let $u,v$ and $x$ be variables such that
\begin{equation}\label{cyc-rep}
K^\ell=u^\ell v^\ell , \ L^\ell=u^{-\ell}v^\ell, \ c=u(x+x^{-1}),
\ E^\ell=y^\ell
\end{equation}
The following representation of $\ue$ in $\C^\ell$ with the linear
basis $\{v_m\}_{m=1}^\ell$ is irreducible:
\begin{eqnarray}\label{IRR}
K v_m&=&uv\e^{2m}v_m, \ Lv_m=u^{-1}v \e^{2m}v_m, \ Ev_m=yv_{m+1}\\
Fv_m&=&y^{-1}u(xv^{-1}\e^{-2m+1}-1)(v\e^{2m-1}-x^{-1})v_{m-1}
\end{eqnarray}
Indeed, since $K$ and $L$ are diagonal in this representation,
they are diagonal in any subrepresentation. Taking into account
the fact that $E$ and $F$ act as cyclic shift operators in the eigen-basis
of $K$ and $L$, one deduces that the only possible subrepresentaions
are the trivial one and the representation itself.
It is easy to see that in this representation
\[
F^\ell=u^\ell(x^\ell v^{-\ell}-1)(v^\ell-x^{-\ell})
\]
The parameters (\ref{cyc-rep}) of the irreducible representation  are
related to the coordinates $a,b,c,d$ on $GL_2^*$ by the following formulae
\[
a=u^\ell v^\ell, d=u^{-\ell}v^\ell, c=-y^\ell u^{-\ell}v^\ell,
b=y^{-\ell}u^\ell(x^\ell-v^\ell-v^{-\ell}+x^{-\ell})
\]
\begin{rem} We have the identities
\[
Tr(b_+(b_-)^{-1})=a+d^{-1}-bcd^{-1}=u^\ell(x^\ell-x^{-\ell})
\]
and
\[
det(b_+(b_-)^{-1})=ad^{-1}=u^{2\ell}
\]
\end{rem}
\begin{prop}\begin{enumerate}
\item Representations with different branches of $v\to v^\ell$ and
$y\to y^\ell$ are isomorphic.

\item Branches of $u\to u^\ell$, $x\to x^\ell$ parameterize
isomorphism classes of irreducible representations for generic
values of $a,b,c,d$.
\end{enumerate}
\end{prop}
The proof is straightforward.

Let $A$ and $B$ be $\ell\times \ell$ matrices acting on the
standard basis in $\C^\ell$ as follows:
\[
Av_n=\e^{2n}v_n, \ Bv_n=v_{n+1}
\]
which satisfy the relations
\[
A^\ell=1, \ B^\ell=1, \ AB=\e^2 BA,
\]
Then the representation (\ref{IRR}) can be written as the
following homomorphism of algebras
\begin{eqnarray}
K&\mapsto& uvA, L\mapsto u^{-1}vA, E\mapsto yB,\\
F&\mapsto& y^{-1}u(xv^{-1}\e A^{-1}-1)(vA\e^{-1}-x^{-1})B^{-1}
\end{eqnarray}
The collection of these homomorphisms is a trivialization
$\phi^i_x: A^i_x\to Math(\ell)$ over a sufficiently small neighborhood
of 1 in $GL_2^*$.

\subsection{The braiding for generic irreducible modules}

The action of $\rc$ on the tensor product of two irreducible
representations is the specialization of the formulae from
section~\ref{raut}. Substituting this action into (\ref{inbr}) we
arrive to the following equations for $R$
\begin{multline}\label{R}
R(1\oti A)R^{-1}=\uto\vto\tuto^{-1}\tvto^{-1}(1\oti A) \\ 
\times(1-\e
\uo^{-1}\vo^{-1}y_1y_2^{-1}\tvto^{-1}A^{-1}B \oti (x_2\vto^{-1}\e
A^{-1}-1)(\vto A\e^{-1}-x_2)B^{-1}A))^{-1},
\end{multline}
\begin{multline}
\tilde{y}_2^{-1}\tuto R(1\oti(\tilde{x}_2\tvto^{-1}\e
A^{-1}-1)(\tvto A\e^{-1}-\tilde{x}_2)B^{-1})R^{-1}\\
=
\uo^{-1}\vo^{-1}(A^{-1}\oti y_2^{-1}\uto(x_2\vto^{-1}\e
A^{-1}-1)(\vto A\e^{-1}-x_2)B^{-1}),
\end{multline}
\begin{equation}
\tuo\tvo\tuto\tvto  R(A\oti A)R^{-1}=\uo\vo\uto\vto (A\oti A),
\end{equation}
\begin{equation}\label{RR}
\tilde{y}_1R(B\oti 1)R^{-1}=y_1\uto\vto(B\oti A),
\end{equation}
Here $\tuo^2=\uo^2$ and $\tuto^2=\uto^2$.

For sufficiently small $x$, $v$, and $z$ such that $z^\ell=(x^\ell
v^{-\ell}-1)(v^\ell-x^\ell)$ define a linear operator
$U:\C^\ell\to \C^\ell$
\[
Uv_n=z\prod_{m=1}^n
((xv^{-1}\e^{-2m+1}-1)^{-1}(v\e^{2m-1}-x^{-1})^{-1}v_n
\]
Let $P_n$ be the projector on the vector $v_n$ in $\C^\ell$:
\[
P_nv_m=\delta_{n,m}v_m,
\]
Define an operator $D$  acting in $\C^\ell\oti\C^\ell$ as
\[
D(v_n\oti v_m)=\e^{2nm}\chi_1^{-n}\chi_2^{m}(v_n\oti v_m)
\]
Assume that the parameters $\chi_1, \chi_2, \e^a$ satisfy the
following relations
\[
u_1v_1u_2v_2=\chi_1\e^{-2a}\tuo\tvo\tuto\tvto
\]
\[
y_1u_2v_2\chi_2=\tilde{y}_1
\]
\[
u_1^{-1}v_1^{-1}z_2\chi_2\e^{-2a}=\tilde{z}_2\tilde{y}_2^{-1}\tuto
\]
and define a power series in $z$ given by the expansion of the function
\begin{equation}\label{PHI}
\Phi(z)=\prod_{m=1}^\ell (1-\e^{2m} z)^{-\frac{m}{\ell}}
\end{equation}
\begin{thm} For generic $x,y\in GL_2^*$ any solution to the equations
\eqref{R}---\eqref{RR} is a scalar multiple of
\begin{equation}\label{finR}
R=D(B^a\oti U)R_1(1\oti {\tilde{U}}^{-1})
\end{equation}
where $R_1=\Phi(sB\otimes B^{-1}\e^{-2})$,
$s=u_2v_2\tuto^{-1}\tvto^{-1}$, $D$ is defined above and $U$ and
$\tilde{U}$ depend on parameters $z_2,u_2, v_2$, and
$\tilde{z}_2,\tilde{u}_2,\tilde{u}_1$, respectively.
\end{thm}
\begin{proof}
The operator $U$ satisfies the identity
\[
U^{-1}(xv^{-1}\e A^{-1}-1)(v A\e^{-1}-x)B^{-1}U=B^{-1}
\]
Write $R$ as
\[
R=D(B^a\oti U)R_1(1\oti {\tilde{U}}^{-1})
\]
Here operator $\tilde{U}$ depends on $\tilde{z}_2, \tuto,
\tvto$ and $U$ depends on $z_2, \uto,\vto$.
Then, the equations for $R$ imply the following relations
for  the operator $R_1$:
\begin{equation}\label{I1}
R_1(A\oti A)R_1^{-1}=A\oti A \ ,
\end{equation}
\begin{equation}\label{I2}
R_1(1\oti A)R_1^{-1}=t(1\oti A)(1-s B\oti B)^{-1}
\end{equation}
\begin{equation}\label{I3}
R_1(1\oti B^{-1})R_1^{-1}=1\oti B^{-1}
\end{equation}
\begin{equation}\label{I4}
R_1(B\oti 1)R_1^{-1}=B\oti 1,
\end{equation}
where
\[
s=u_2v_2\tuto^{-1}\tvto^{-1}, \
t=\e\chi_2\chi_3u_1^{-1}v_1^{-1}y_1y_2^{-1}\tvto^{-1}
\]
Identities \eqref{COP}---\eqref{COPP} imply that $t^\ell=1-s^\ell$.

The function $\Phi$ satisfies the following "difference equation":
\begin{equation}\label{deq}
\Phi(\e^2 z)=(1-z^\ell)^{\frac{1}{\ell}}(1-\e^2z^{-1})^{-1}\Phi(z)
\end{equation}
From equations (\ref{I1}), (\ref{I3}), (\ref{I4}) we conclude that
$R_1$ is a polynomial in $B\otimes B^{-1}$. Equation~(\ref{I2}) 
determines the coefficients of this polynomial up to an
overall scalar factor. Taking into account the difference equation
for $\phi$, we conclude that there is one dimensional family of
solutions to linear equations~(\ref{I2}):
\[
R_1=const \ \Phi(sB\otimes B^{-1}\e^{-2})
\]
\end{proof}
\begin{rem}
The power series $\Phi(z)$ can be analytically continued to the meromorphic
function $\Phi(z)$ defined in a neighborhood of $z=0$. More
precisely $\Phi$ determines the  curve $C_\Phi=\{
\Phi^\ell\prod_{m=1}^\ell (1-\e^{2m} z)^{-\frac{m}{\ell}}=1 \}$ in
$\C\times \C$. The difference equation is the action of $\Z_\ell$
on the curve $C_1=\{ (\Phi, \e,\phi)\in \C^3|
\Phi^\ell\prod_{m=1}^\ell (1-\e^{2m} z)^{-\frac{m}{\ell}}=1,
\phi^\ell+z^\ell=1\}$. The element $\omega\in \Z_\ell$ acts as
$(\Phi, \e,\phi)\to (\Phi\phi(1-\omega z)^{-\ell},z\omega, \phi)$.
\end{rem}
It is not difficult to compute the determinant of the matrix $R$:
\[
det(R)=(1-s^\ell)^{\frac{\ell(\ell+2)}{2}} \ ,
\]
Since $R$ is fixed yet only up a scalar function
of parameters of the representations, then the element
\[
\tilde{R}=(1-s^\ell)^{\frac{\ell+1}{2\ell}}R \,
\]
intertwines the same representations as $R$ with the property
$det(\tilde{R})=1$ and at $s=0$ the matrix
$\tilde{R}=D_{12}B_1^aU_2\tilde{U}_2^{-1}$ satisfies the
Yang--Baxter equation. Therefore the factor $c^{ijk}$ is trivial
for the matrix $\tilde{R}$.

Thus, every pair of generic
irreducible $\ue$-modules defines the following solution to the holonomy
Yang--Baxter equation:
\[
R^{i,j}(x,y)=(1-s^\ell)^{\frac{\ell+1}{2\ell}}D(B^a\otimes
U)\Phi(sB\otimes B^{-1}\e^2)(1\otimes \tilde{U}^{-1}) \ .
\]
Here $x$ and $y$ are in sufficiently small neighborhood of 1 in
$GL_2^*$ and all other ingredients are defined above.

\appendix

\section{}\label{double}

{\bf 1.} Let $b \subset sl_2$ be the Borel subalgebra and $H$,
the generator of the Cartan subalgebra.  The algebra $U_hb$ is
generated by $H$ and $X$ with the defining relation
\[
HX - XH = 2X.
\]
This is a Hopf algebra over ${\mathbb C}[[h]]$ with the
comultiplication
\begin{eqnarray*}
\D H &= &H \otimes 1 + 1 \otimes H,\\
\D X &= &X \otimes e^{\frac {hH}{2}} + 1 \otimes X.
\end{eqnarray*}

{\bf 2.} Let $(U_hb)^0$ be the Hopf algebra generated by
elements $H^\vee,\ X^\vee  $, completed by formal power series in
$h$ and $H^\vee$ with the defining relations
\begin{eqnarray*}
[H^\vee  ,X^\vee  ] &= &-\frac {h}{2} X^\vee  ,\\
\D H^\vee   &= &H^\vee   \otimes 1 + 1 \otimes H^\vee  ,\\
\D X^\vee   &= &X^\vee   \otimes 1 + e^{-2H^\vee  } \otimes X^\vee  ,
\end{eqnarray*}

{\bf 3.} There is a nondegenerate pairing $\langle
\cdot,\cdot\rangle\colon U_hb \otimes (U_hb)^0 \rightarrow {\mathbb
C}[[h]]$ of Hopf algebras with
\begin{eqnarray}\label{mul}
\langle l,ba\rangle &= &\langle \Delta l,a \otimes b\rangle \\
\langle lm,a\rangle &= &\langle l \otimes m,\D a\rangle
\end{eqnarray}
This pairing is defined on generators as
\begin{eqnarray*}
\langle H^\vee  ,H\rangle &= &1,\ \langle X^\vee  ,X\rangle = 1\\
\langle H^\vee  ,X\rangle &= &\langle X^\vee  ,H\rangle = 0.
\end{eqnarray*}
It is easy to extend this pairing to the whole algebra using
(\ref{mul}). The pairing between the monomials is
\[
\langle (H^\vee  )^n(X^\vee  )^m,H^{n'},X^{m'}\rangle = 
\delta_{nn'}\delta_{mm'}n!b_m
\]
where
\[
b_m = \frac {1-e^{-hn}}{1-e^{-h}} \cdot \frac {1-e^{-h(n-1)}}{1-e^{-h}}
\dots \frac {1-e^{-h}}{1-e^{-h}}.
\]
Indeed, it is clear that the right-hand side is proportional to
$\delta_{n,n'}\delta_{mm'}$.  Let us commute the coefficient.
Using the Hopf properties (\ref{mul}) of the pairing and the
grading arguments we obtain:
\begin{eqnarray*}
\langle (H^\vee  )^n(X^\vee  )^m, H^nX^m\rangle &= &\langle 
(H^\vee  )^n \otimes
(X^\vee  )^m,(\D H)^n(\D X)^m\rangle \\
&= &\langle (H^\vee  )^n \otimes (X^\vee  )^m,H^n \otimes X^m\rangle \\
&= &\langle (H^\vee  )^n,H^n\rangle\langle (X^\vee  )^m,X^m\rangle.
\end{eqnarray*}
For the first factor we have
\begin{eqnarray*}
\langle (H^\vee  )^n,H^n\rangle &= &\langle (H^\vee  )^{n-1} 
\otimes H^\vee  ,(\D
H)^n\rangle \\
&= &\langle (H^\vee  )^{n-1} \otimes H^\vee  ,nH^{n-1} \otimes H\rangle \\
&= &n\langle (H^\vee  )^{n-1},H^n\rangle
\end{eqnarray*}
therefore
\[
\langle (H^\vee  )^n,H^n\rangle = n!
\]
For the second factor we have
\begin{eqnarray*}
\langle (X^\vee  )^n,X^n\rangle &= &\langle (X^\vee  )^{n-1} 
\otimes X^\vee  ,(\D
X)^n\rangle \\
&= &\langle (X^\vee  )^{n-1} \otimes X^\vee  , 
\frac {-e^{-hn}}{1-e^{-h}} X^{n-1}
\otimes e^{\frac {hH}{2}(n-1)}X\rangle.
\end{eqnarray*}
Here we used the second term in
\begin{eqnarray*}
(\D X)^n &= &(X \otimes e^{\frac {hH}{2}} + 1 \otimes X)^n \\
&= &X^n\otimes e^{\frac {hnH}{2}} \\
&+ &(X^{n-1} \otimes e^{\frac {h(n-1)H}{2}}X + X^{n-1} \otimes e^{\frac
{hH}{2}} X e^{\frac {h(n-2)}{2}H} + \dots + X^{n-1} \otimes X e^{\frac
{h(n-1)H}{2}}) + \dots \\
&= &X^n \otimes e^{\frac {hnH}{2}} + (1 + e^{-h} + \dots + e^{-(n-1)h})
X^{n-1} \otimes e^{\frac {(n-1)}{2}hH} X + \dots
\end{eqnarray*}
Therefore
\[
\langle (X^\vee  )^n,X^n\rangle = \frac {1-e^{-hn}}{1-e^{-h}} \langle
(X^\vee  )^{n-1},X^{n-1}\rangle \cdot \langle X^\vee  ,e^{\frac {hH(n-1)}{2}}
X\rangle.
\]
Taking into account that
\[
\langle X^\vee  ,e^{\frac {hH}{2}(n-1)}X\rangle 
= \langle \D X^\vee  ,X \otimes
e^{\frac {hH}{2}(n-1)}\rangle = 1
\]
we obtain
\[
\langle (X^\vee  )^n,X^n\rangle = b_n.
\]
where $b_n$ is defined above.

{\bf 4.} The double of the Hopf algebra $U_hb$ is the Hopf
algebra structure on the space ${\mathcal
D}(U_hb) = U_hb \otimes (U_hb)^0$ completed in formal power series
in $h$ and $H^\vee$ such that the coalgebra
structure is the tensor product of coalgebras.  The algebra
structure is completely determined by the condition that the natural
embeddings $U_hb,(U_hb)^0 \hookrightarrow {\mathcal D}(U_hb)$ are
Hopf algebra homomorphisms and that the canonical element $R =
\sum_i e_i \otimes e^i$ intertwines the coproduct with the
opposite coproduct
$$
R\D(a) = \D^{op}(a)R,
$$
where $\D^{op}(a) = \s_0\D(a)$, $\s(x \otimes y) = y \otimes x$.

Equivalently, the product in ${\mathcal D}(U_hb)$ can be defined as

\begin{equation}
(a\otimes l)(b\otimes m)=\sum ab_{(2)}\otimes l_{(2)}m
<b_{(1)}, S^{-1}_{A^*}(l_{(1)})><b_{(3)}, l_{(3)}>
\end{equation}
Here we use the notation $\D(a)=\sum a_{(1)}\otimes a_{(2)}$
for the comultiplication in $U_hb$
and $\D^o(l)=\sum l_{(2)}\otimes l_{(1)}$ for the comultiplication
in $(U_hb)^o$ (see \cite{Ma}).

Using either of these definitions, it is easy to show that ${\mathcal
D}(U_hb)$ is isomorphic to the algebra generated by $H,X,H^\vee
,X^\vee  $, completed in formal power series in $h$ and $H^\vee$
with the defining relations
\begin{eqnarray*}
\ [H,H^\vee  ] &= &0,\ [H,X] = 2X,\ [H^\vee  ,X] = \frac {h}{2} X \\
\ [H,X^\vee  ] &= &-2X^\vee  ,\ [H^\vee  ,X^\vee  ] 
= -\frac {h}{2} X^\vee   \\
\ [X,X^\vee  ] &= &e^{\frac {hH}{2}} - e^{-2H^\vee  }.
\end{eqnarray*}
The universal $R$-matrix is the canonical element $R = \sum_i e_i \otimes
e^i \in U_hb \otimes (U_hb)^0$
\begin{eqnarray}\label{Rdoubl}
R &= &\sum_{n,m \ge 0} \frac {1}{n!b_m} H^nX^m 
\otimes (H^\vee  )^n(X^\vee  )^m \\
&= &\exp(H \otimes H^\vee  )f(X \otimes X^\vee  ;e^{-h})
\end{eqnarray}
where
\[
f(z;q) = \sum_{n \ge 0} \frac {(1-q)^nz^n}{(1-q)\dots(1-q^n)} =
\prod_{n \ge 0} (1-(1-q)zq^n)^{-1}.
\]
Here $b_m$ is as above and $q=e^{-h}$. It follows immediately from
the definition of the double that
$$
\begin{array}{rll}\label{DR}
(\D \otimes id)(R) &= &R_{13}R_{23} \\
(id \otimes \D)(R) &= &R_{13}R_{12}
\end{array}
$$
and, in particular, that $R$ satisfies the Yang--Baxter equation.

{\bf 5.} The quantum $gl_2$, or $U_hgl_2$ is the Hopf algebra
generated by $H,G,X,Y$ with the defining relations and the
coproduct described in (\ref{u1})--(\ref{u3}).  It is clear that
the mapping $\varphi: {\mathcal D}(U_hb) \rightarrow U_hgl_2$
acting on generators as
\[
\varphi(H) = H,\ \varphi(H^\vee  ) = \frac {h}{4} G,\ \varphi(X) = X,\
\varphi(X^\vee ) = Y(e^{\frac{h}{2}}-e^{-\frac{h}{2}})
\]
extends to a Hopf algebra homomorphism.  The element (\ref{Rh} )
is the image of the universal $R$-matrix (\ref{Rdoubl}) and therefore
satisfies (\ref{br}).

\section{}\label{INT}

In section \ref{raut} we defined the action of $\rc$ on the
generators $E\otimes 1$, $1\otimes F$, $1\oti K$, and $1\oti L$
explicitly and on other generators implicitly by requiring the intertwining
property (\ref{rcom}).
The action of $\rc$ on element $F\oti 1$ and on $1\oti E$ can also
be found from the condition that $\rc(c_1\oti c_2)= c_1\oti c_2$
for central elements $c_1$ and $c_2$ and $\rc(K\oti K)=K\oti K$.
Let us verify this.

First, notice that $\rc(K\oti K)=K\oti K$ together with
(\ref{rc-1}) implies that
\[
\rc(K\oti 1)=(1-tK^{-1}E\oti FL)(K\oti1)
\]
Now let us find the action of $\rc$ on $1\oti E$:
\begin{multline}
\rc(1\oti E)=\rc(1\oti EF)(K\oti F^{-1})=\\
\rc(1\oti(c-t^{-1}K-tL^{-1}))(K\oti F^{-1})=\\
1\oti cF^{-1}-(t^{-1}\rc(1\oti K)+t\rc(1\oti L^{-1}))(K\oti F^{-1})=\\
1\oti E+t^{-1}K\oti KF^{-1}+tK\oti L^{-1}F^{-1}-\\
t^{-1}(1\oti K)(1-tK^{-1}E\oti FL)^{-1}(K\oti
F^{-1}-t(1-tK^{-1}E\oti FL) (1\oti L^{-1}F^{-1})\\=1\oti E+E\oti
K-(1-t^{-1}K^{-1}E\oti FL)^{-1}E\oti KL,
\end{multline}
Similarly, we have:
\begin{multline}
\rc(F\oti 1)=(E^{-1}\oti L^{-1})\rc(EF\oti 1)=\\
(E^{-1}\oti L^{-1})\rc((c-t^{-1}K-tL^{-1})\oti 1)=\\
E^{-1}c\oti L^{-1}-t^{-1}(E^{-1}\oti L^{-1})(\rc(K\oti 1)
+t\rc(L^{-1}\oti 1))=\\
F\oti L^{-1}+t^{-1}E^{-1}K\oti L^{-1} +tE^{-1}L^{-1}\oti L^{-1}-\\
t^{-1}(E^{-1}\oti L^{-1})(1-tK^{-1}E\oti FL)(K\oti
1)-t(E^{-1}L^{-1}\oti L^{-1})(1-tK^{-1}E\oti FL)\\=F\oti
L^{-1}+1\oti F-(KL^{-1}\oti F) (1-tK^{-1}E\oti FL)^{-1}.
\end{multline}
Taking into account the action of $\rc$ on $K\oti F$ and $E\oti
L^{-1}$, it is easy to see that the previous formulae are
equivalent to
\[
\rc(E\oti K+1\oti E)=K\oti E+E\oti 1,
\]
\[
\rc(F\oti 1+L^{-1}\oti F)=1\oti F+ F\oti L^{-1}
\]

\end{document}